\documentclass{amsart}

\usepackage{amsmath,amssymb,amsthm}
\usepackage{graphicx,tikz}

\newtheorem*{thm}{Theorem}

\theoremstyle{definition}

\theoremstyle{remark}

\begin{document}

\title[]{Varadhan asymptotics for the \\heat kernel on finite graphs}

\author[]{Stefan Steinerberger}
\address{Department of Mathematics, Yale University, New Haven, CT 06511, USA}
\email{stefan.steinerberger@yale.edu}

\begin{abstract} Let $G$ be a simple, finite graph and let $p_t(x,y)$ denote the heat kernel on $G$. The purpose of this short note is to show that
for $t \rightarrow 0^+$
$$ p_t(x,y) = \# \left\{\mbox{paths of length}~d(x,y)~\mbox{between}~x~\mbox{and}~y\right\} \frac{t^{d(x,y)}}{d(x,y)!} + \mathcal{O}(t^{d(x,y)+1}),$$
where $d(x,y)$ is the usual Graph distance.
This is the discrete analogue of the classical Varadhan asymptotic for the heat kernel on manifolds and refines a result of Keller, Lenz, M\"unch, Schmidt \& Telcs. 
The asymptotic behavior encapsulates additional geometric information: if the Graph is bipartite, then the next term in
the expansion is negative.
\end{abstract}

\maketitle
\section{Introduction and main result}
\subsection{Introduction.} The geometry of a Riemannian manifold $(M,g)$ is strongly coupled with the behavior of the diffusion process induced by the
Laplace-Beltrami operator $\Delta_g/2$. One of the foundational examples of this connection is a seminal result of Varadhan
showing that the behavior of the
$$ \mbox{heat kernel on}~\mathbb{R}^n \qquad p_t(x,y) = \frac{1}{(2 \pi t)^{n/2}} e^{-\frac{\|x-y\|^2}{2t}}$$
remains valid for short time on manifolds in the sense of
$$ \lim_{t \rightarrow 0^+}{ t \log{p_t(x,y)}} = - \frac{d_g(x,y)^2}{2}.$$
The behavior of the heat kernel on graphs is a rather large subject, we refer to 
Chung \cite{chung}, Davies \cite{davies}, Metzger \& Stollmann \cite{metz}, Pang \cite{pang} and Schmidt \cite{schm}.
The natural question of short-time asymptotics has only been addressed recently by
Keller, Lenz, M\"unch, Schmidt \& Telcs \cite{keller}. They show that, as $t \rightarrow 0^+$, the behavior
of the heat kernel is given by
$$ | p_t(x,y) - c(x,y) t^{d(x,y)} | \leq C(x,y) t^{d(x,y)+1},$$
where $c(x,y)$ and $C(x,y)$ are two numerical constants depending on $x,y$. This, in particular, implies
$$ \lim_{t \rightarrow 0^+} \frac{ \log{ p_t(x,y)}}{ \log{t}} = d(x,y)$$
where $d(x,y)$ is the usual combinatorial graph distance. As a consequence, the behavior of the heat kernel
on a graph for short times is substantially different from the behavior on a manifold (while, for longer time, one would expect
them to behave in a similar manner). The proof in \cite{keller} is based on an abstract result in functional
analysis and also applies to certain infinite (locally finite) graphs.

\subsection{The result.} The purpose of this paper is to give a short and simple proof of this result for finite, simple graphs that also determines
the leading coefficient in the asymptotic expansion. This coefficient turns out to be
$$ c(x,y) = \frac{1}{d(x,y)!}\# \left\{ \mbox{paths of length}~d(x,y)~\mbox{between}~x~\mbox{and}~y\right\}$$
yet another natural geometric quantity that further emphasizes the interplay between geometry and
analysis. Let $G=(V,E)$ be a finite graph of cardinality $|V|=n$ and identify vectors
in $\mathbb{R}^n$ with functions $f:V \rightarrow \mathbb{R}$. The Graph Laplacian $\mathcal{L}$ is 
a linear operator acting on functions via
$$ (\mathcal{L}f)(x) = \sum_{(x,y) \in E}{(f(y) - f(x))}.$$
The linear operator $\mathcal{L}$ can be identified with the so-called Kirchhoff matrix and introduces a natural diffusion $e^{t \mathcal{L}}f$ defined as solution of the heat equation
$
(\partial_t - \mathcal{L})e^{t\mathcal{L}}f = 0.$
The heat kernel is defined as 
$ p_t(x,y) = \left( e^{t \mathcal{L}} \delta_x\right)(y).$
We observe $p_t(x,y) = p_t(y,x)$ and, by linearity,
$$ \left( e^{t \mathcal{L}}f \right)(x) = \sum_{y \in V}{p_t(x,y) f(y)}.$$

\begin{thm} Let $G$ be a finite, simple graph. Then, as $t \rightarrow 0^+$, we have
$$ p_t(x,y) = \# \left\{ \mbox{paths of length}~d(x,y)~\mbox{between}~x~\mbox{and}~y \right\} \frac{t^{d(x,y)}}{d(x,y)!} + \mathcal{O}(t^{d(x,y)+1}).$$
 If $G$ is bipartite, then the constant in front of
the $t^{d(x,y)+1}-$term is negative. 
\end{thm}

We quickly illustrate the main result on a simple example (see Fig. 1): there are two paths of length
$d(a,b) = 2$ going from $a$ to $b$ and three paths of length $d(a,c) = 3$ going from $a$ to $c$.
An explicit computation shows that for $t$ small
$$ p_t(a,b) = 2\frac{t^2}{2} - \frac{5}{2}t^3 + \mathcal{O}(t^4) \quad \mbox{and} \quad p_t(a,c) = 3\frac{t^3}{6} - \frac{7}{6}t^4 + \mathcal{O}(t^5).$$

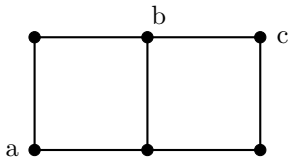
\begin{figure}[h!]
\begin{center}
\begin{tikzpicture}[scale=1.5]
\filldraw (0,0) circle (0.05cm);
\filldraw (1,0) circle (0.05cm);
\filldraw (2,0) circle (0.05cm);
\filldraw (0,1) circle (0.05cm);
\filldraw (1,1) circle (0.05cm);
\filldraw (2,1) circle (0.05cm);
\draw [thick] (0,0) -- (2,0);
\draw [thick] (0,1) -- (2,1);
\draw [thick] (0,0) -- (0,1);
\draw [thick] (1,0) -- (1,1);
\draw [thick] (2,0) -- (2,1);
\node at (-0.2, 0) {a};
\node at (1.1, 1.2) {b};
\node at (2.2, 1) {c};
\end{tikzpicture}
\end{center}
\caption{A bipartite graph with three labeled vertices.}
\end{figure}

The result could be phrased in a different language. Let $|V|=n$, let $\lambda_1, \dots, \lambda_n$ denote the eigenvalues
of $\mathcal{L}$ and $v_1, \dots, v_n$ denote the $\ell^2-$normalized eigenvectors. $\delta_x$ denotes a vector that is 1 in the $x-$th coordinate and
0 everywhere else. Then
$$ p_t(x,y) = \sum_{k=1}^{n}{ e^{-\lambda_k t} \left\langle v_k, \delta_x\right\rangle \left\langle v_k, \delta_y\right\rangle}.$$
Our main result can then be formulated as the identity
$$ \sum_{k=1}^{n}{ (-\lambda_k)^{d(x,y)} \left\langle v_k, \delta_x\right\rangle \left\langle v_k, \delta_y\right\rangle} = \# \left\{ \mbox{paths of length}~d(x,y)~\mbox{between}~x~\mbox{and}~y \right\}.$$

\subsection{Remarks.}  As shown by Varadhan \cite{vara1,vara2}, Brownian motion
in the small-time limit is moving close to geodesics whenever $x$ and $y$ have a unique geodesic joining them. Our observation mirrors the very same phenomenon (with the difference that
most points are connected by more than one geodesic): the proof shows that as $t \rightarrow 0$, all of $p_t(x,y)$ is flowing
exactly along the paths of minimal length from $x$ to $y$. This shows the case of finite graphs could even be regarded as a basic toy example motivating Varadhan's result in the continuous setting. 

 \begin{center}
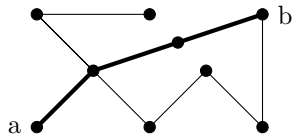
\begin{figure}[h!]
\begin{tikzpicture}[scale=1.5]
\filldraw (0,0) circle (0.05cm);
\filldraw (1,0) circle (0.05cm);
\filldraw (0,1) circle (0.05cm);
\filldraw (1,1) circle (0.05cm);
\filldraw (-1,0) circle (0.05cm);
\filldraw (-1,1) circle (0.05cm);
\filldraw (-0.5,0.5) circle (0.05cm);
\filldraw (0.5,0.5) circle (0.05cm);
\filldraw (0.25,0.75) circle (0.05cm);
\draw [] (-1, 0) -- (-0.5, 0.5) -- (-1, 1) -- (0,0);
\draw [] (-0.5, 0.5) -- (1,1);
\draw [] (-1, 1) -- (0,1);
\draw [] (-0.5, 0.5) -- (1,1);
\draw [] (0, 0) -- (0.5,0.5) -- (1,0) -- (1,1);
\draw [ultra thick] (-1,0) -- (-0.5, 0.5) -- (1,1);
\node at (-1.2, 0) {a};
\node at (1.2, 1) {b};
\end{tikzpicture}
\caption{A practical application of the result: to find the graph distance $d(a,b)$ it suffices to compute $p_t(a,b)$ for small $t$. The coefficient determines the number of paths with minimal length. (This idea is due to Crane, Weischedl \& Wardetzky \cite{cr1,cr2}.)}
\end{figure}
\end{center}
\vspace{-10pt}

Finally, we note that the result has some immediate applications. Crane, Weischedl \& Wardetzky \cite{cr1,cr2} propose solving
the heat equation and using Varadhan's asymptotic expansion as a way of computing geodesic distances. Our formula provides a theoretical
justification for this method when discretizing a manifold with a mesh; it remains to be seen whether it can helpful in parameter selection or mesh design.
  
 \section{Proof}
 \begin{proof} The proof proceeds in two steps; we first establish the asymptotic behavior and then use this
 information to re-run the argument and obtain more refined information in the case where $G$
 is bipartite. We fix $x \in V$. The main idea is to reinterpret the problem as a coupled system of ordinary differential equations. 
 For this linear system, we have an explicit solution given by the matrix exponential applied to the Kirchhoff matrix $L$ and
 $ e^{t \mathcal{L}} \delta_x  = e^{t L} \delta_x$.
 The main ingredient, similar to \cite{keller}, is Taylor's theorem which we use in the form
  $$ e^{tL} = 1 + t L + \mathcal{O}(t^2).$$
 This is already sufficient to prove that for $d(x,y) = 1$ 
 $$ p_t(x,y) = t + \mathcal{O}(t^2)$$
 and that $p_t(x,y) = \mathcal{O}(t^2)$ for all $d(x,y) \geq 2$. We will now use this as basis for a proof
 by induction: the precise statement is the Theorem as desired
 and that $p_t(x,y) = \mathcal{O}(t^{n+1})$ for all $d(x,y) \geq n+1$. We assume this to hold for $d(x,y) \leq n$. Let now $y \in V$ be
 chosen such that $d(x,y) = n+1$. By definition
 $$ \partial_t  \left[ e^{t \mathcal{L}} \delta_x \right] (y)  = \left[ \mathcal{L} e^{t \mathcal{L}} \delta_x \right](y) = \sum_{(y,z) \in E}{ \left[e^{t\mathcal{L}}\delta_x\right](z) - \left[e^{t\mathcal{L}}\delta_x\right](y)}.$$ 
 The fundamental theorem of calculus implies
 $$  \left[ e^{t \mathcal{L}} \delta_x \right] (y) = \int_{0}^{t}{  \sum_{(y,z) \in E}{ \left[e^{s\mathcal{L}}\delta_x\right](z) - \left[e^{s\mathcal{L}}\delta_x\right](y)} ds}.$$
 By induction assumption, we have $ \left[e^{s\mathcal{L}}\delta_x\right](y) = \mathcal{O}( s^{n+1})$ and integrating thus turns the second term into one of size $t^{n+2}$ which is of lower order. The
 same reasoning can be applied to neighbors of $z$ of $y$ that are not closer to $x$ than $y$, i.e. $(y,z) \in E$, satisfying $d(x,z) \geq n+1$. This yields 
  $$  \left[ e^{t \mathcal{L}} \delta_x \right] (y) = \int_{0}^{t}{  \sum_{(y,z) \in E \atop d(x,z) = n}{ \left[e^{s\mathcal{L}}\delta_x\right](z)ds} } + \mathcal{O}(t^{n+2}).$$
We can now invoke the induction hypothesis and obtain
 $$  \left[ e^{t \mathcal{L}} \delta_x \right] (y) = \frac{t^{n+1}}{(n+1)!} \sum_{(y,z) \in E \atop d(x,z) = n}{  (\# \mbox{paths of length}~n~\mbox{from}~x~\mbox{to}~z)} + \mathcal{O}(t^{n+2}).$$ 
However, the sum simplifies to the total number of paths of length $n+1$ from $x$ to $y$. It remains to show that $ \left[e^{t\mathcal{L}}\delta_x\right](y) = \mathcal{O}( t^{n+2})$ for all vertices $y$ satisfying $d(x,y) \geq n+2$ to complete the induction step (we already know that they are $\mathcal{O}(t^{n+1})$). This follows easily from another application of
 $$  \left[ e^{t \mathcal{L}} \delta_x \right] (y) = \int_{0}^{t}{  \sum_{(y,z) \in E}{ \left[e^{s\mathcal{L}}\delta_x\right](z) - \left[e^{s\mathcal{L}}\delta_x\right](y)} ds}$$
and the fact that the all neighbors $z$ of $y$ satisfy $d(x,z) \geq n+1$.

The second part of the statement is again shown by induction. We have already determined the leading behavior,
the additional statement is that the coefficient in front of the next term of the expansion is negative and we prove this with induction
over $n$ for all vertices $y$ satisfying $d(x,y) \leq n$ (where, as before, $x$ is fixed without loss of generality). The case $n = 0$
is simple since
$$   \left[ e^{t \mathcal{L}} \delta_x \right] (x) = 1 - \mbox{deg}(x)t + \mathcal{O}(t^2).$$ 
The case $n=1$ is the first interesting case. We use Taylor's theorem
 $$ e^{tL} = 1 + t L + \frac{t^2 L^2}{2} + \mathcal{O}(t^3).$$
Let $d(x,y) = 1$. Since $G$ is bipartite, the vertices $x$ and $y$ have no common neighbors and the relevant entry of the square of the Kirchhoff matrix
is easy compute 
$$   \left[ e^{t \mathcal{L}} \delta_x \right] (y) = t  -  (\mbox{deg}(x) + \mbox{deg}(y))\frac{t^2}{2} +  \mathcal{O}(t^3).$$ 
We now assume the desired statement hold for all vertices at distance at most $n$ from $x$ and assume $d(x,y) =  n+1$. We again
use
\begin{align*}
  \left[ e^{t \mathcal{L}} \delta_x \right] (y) &= \int_{0}^{t}{  \sum_{(y,z) \in E}{ \left[e^{s\mathcal{L}}\delta_x\right](z) - \left[e^{s\mathcal{L}}\delta_x\right](y)} ds}\\
  &= \int_{0}^{t}{  \sum_{(y,z) \in E}{ \left[e^{s\mathcal{L}}\delta_x\right](z) ds}} - \mbox{deg}(y) \int_{0}^{t}{ \left[e^{s\mathcal{L}}\delta_x\right](y) ds}.
  \end{align*}
The second term is easy to deal with: our main result implies
 $$ - \mbox{deg}(y) \int_{0}^{t}{ \left[e^{s\mathcal{L}}\delta_x\right](y) ds} = - \mbox{deg}(y) c \frac{t^{n+2}}{n+2} + \mathcal{O}(t^{n+3}),$$
 where $c>0$ is the leading term in the asymptotic expansion. In particular, this term contributes a negative coefficient at order $n+2$. We shall
 now show that the first term also only contributes negative terms and this then implies the desired statement. We observe that $(y,z) \in E$
 implies that $d(x,z) \in \left\{n, n+2\right\}$ since $G$ is bipartite. By induction assumption
 $$  \int_{0}^{t}{  \sum_{(y,z) \in E \atop d(x,z) = n}{ \left[e^{s\mathcal{L}}\delta_x\right](z) ds}}$$
 contributes a negative coefficient at order $n+2$. As a consequence of our main result
  $$ \int_{0}^{t}{  \sum_{(y,z) \in E \atop d(x,z) = n+2}{ \left[e^{s\mathcal{L}}\delta_x\right](z) ds}} = \mathcal{O}(t^{n+3})$$
  and this implies the desired result.
 \end{proof}

 The result can be easily extended to graphs with weighted edges: the final
result will be that every path has a weight given as the product of the weight of all of its edges and the proof is identical. Returning to the basic idea that diffusion on a finite graph can be written as a coupled linear system of ordinary differential equations, we note that there is an explicit solution given by the matrix exponential
$$ e^{t L} = \sum_{k=0}^{\infty}{\frac{L^k}{k!} t^k},$$
where $L$ is the Kirchhoff matrix. This could be used to give an alternative formulation of the proof in a different language as follows. Let us label the vertices of $V$
by $\left\{1,\dots, |V|\right\}$. Then it suffices to show that, for fixed $i \neq j$, the smallest $k$ for which $(L^k)_{ij} \neq 0$ is exactly given by $k = d(i,j)$ and, moreover, that for $k = d(i,j)$ the entry coincides the entry of the $k-$th power of the adjacency matrix $(L^k)_{ij} = (A^k)_{ij}$ (counting the number of paths). This can be shown by induction as we did above.

\end{document}